\documentclass[12pt]{article}

\usepackage[left=2cm,right=2cm,top=4cm,bottom=4cm]{geometry}
\usepackage[T1]{fontenc}
\usepackage[latin1]{inputenc}
\usepackage{lmodern}
\usepackage[nottoc, notlof, notlot]{tocbibind}
\usepackage{amsmath}
\usepackage{amsthm}
\usepackage{amssymb}
\usepackage{ stmaryrd }
\usepackage{graphicx}
\usepackage{tikz} 
\usepackage{dsfont}
\usepackage{xspace}
\usepackage{enumitem}
\usepackage{ upgreek }
\usetikzlibrary{matrix,arrows} 
\usepackage[colorlinks,urlcolor=blue,linkcolor=blue,citecolor=brown]{hyperref}

\theoremstyle{plain}
\newtheorem{thm}{Theorem}

\theoremstyle{remark}

\theoremstyle{plain}
\newtheorem{prop}{Proposition}

\theoremstyle{definition}
\newtheorem*{defn}{Definition}

\theoremstyle{plain}

\theoremstyle{plain}

\theoremstyle{plain}

\theoremstyle{remark}
\newtheorem*{rmq}{Remark}

\begin{document}

\title{About Jarn\'ik's-type relation in higher dimension}
\author{Antoine MARNAT\footnote{supported by the Austrian Science Fund (FWF), Project F5510-N26, which is a part of the Special Research Program ``Quasi-Monte Carlo Methods: Theory and Applications'' and FWF START-project Y-901} \\
 \href{mailto:marnat@math.tugraz.at}{marnat@math.tugraz.at}}
 
 \date{}
\maketitle

\abstract{  Using the Parametric Geometry of Numbers introduced recently by W.M. Schmidt and L. Summerer \cite{SchSu,SchSu2} and results by D. Roy \cite{RoyParam,RoySpec}, we show that German's transference inequalities between the two most classical exponents of uniform Diophantine approximation are optimal. Further, we establish that the $n$ uniform exponents of Diophantine approximation in dimension $n$ are algebraically independent. Thus, no Jarn\'ik's-type relation holds between them.}

 \section{Introduction}
 
Throughout this paper, the integer $n\geq1$ denotes the dimension of the ambient space, $\boldsymbol{\theta}=(\theta_1, \ldots , \theta_n)$ denotes an $n$-tuple of real numbers such that $1,\theta_1, \ldots , \theta_n$ are $\mathbb{Q}$-linearly independent.\\
 
Let $d$ be an integer with $0\leq d \leq n-1$. We define the exponent ${\omega}_{d}(\boldsymbol{\theta})$  (resp. the uniform exponent $\hat{\omega}_{d}(\boldsymbol{\theta})$) as the supremum of the real numbers $\omega$ for which there exist  rational affine subspaces $L \subset \mathbb{R}^{n}$ such that 
 \[ \dim(L)=d \; ,  \quad H(L)\leq H \quad \textrm{ and } \quad H(L)d(\boldsymbol{\theta},L) \leq H^{-\omega}  \]
  for arbitrarily large real numbers $H$  (resp. for every sufficiently large real number $H$). Here $H(L)$ denotes the height of $L$ (see \cite{SchH} for more details), and $d(\boldsymbol{\theta},L)=\min_{P\in L} d(\boldsymbol{\theta},P)$ is the minimal distance between $\boldsymbol{\theta}$ and a point of $L$. \\

These exponents were introduced originally by M. Laurent \cite{MLwd}. They interpolate between the classical exponents $\omega(\boldsymbol{\theta})=\omega_{n-1}(\boldsymbol{\theta})$ and $\lambda(\boldsymbol{\theta})=\omega_0(\boldsymbol{\theta})$ (resp. $\hat{\omega}(\boldsymbol{\theta})=\hat{\omega}_{n-1}(\boldsymbol{\theta})$ and $\hat{\lambda}(\boldsymbol{\theta})=\hat{\omega}_0(\boldsymbol{\theta})$) that were  introduced by A. Khinchin \cite{Khin2,Khin1}, V. Jarn\'ik \cite{JAR} and Y. Bugeaud and M. Laurent \cite{BugLau2,BugLau}.\\
 
 We have the relations
 \[ \omega_{0}(\boldsymbol{\theta}) \leq \omega_{1}(\boldsymbol{\theta}) \leq \cdots \leq \omega_{n-1}(\boldsymbol{\theta}) ,  \]
 \[\hat{\omega}_{0}(\boldsymbol{\theta})\leq \hat{\omega}_{1}(\boldsymbol{\theta}) \leq \cdots \leq \hat{\omega}_{n-1}(\boldsymbol{\theta}),\]
 and Minkowski's First Convex Body Theorem \cite{Mink} and Mahler's compound convex bodies theory provide the lower bounds
 
 \begin{equation*}
 \omega_{d}(\boldsymbol{\theta}) \geq \hat{\omega}_{d}(\boldsymbol{\theta}) \geq \cfrac{d+1}{n-d} , \quad \textrm{  for  } \; 0 \leq d \leq n-1 .
 \end{equation*}
 
 These exponents happen to be related, as was first noticed by Khinchin with his transference theorem \cite{Khin1}. The study of these transferences has two aspects. First, establishing transference inequalities valid for every suitable point $\boldsymbol{\theta}$. Then, there is the reverse problem, that consists in constructing points $\boldsymbol{\theta}$ to show that these inequalities are sharp. For this, one can prove that there exists points $\boldsymbol{\theta}$ whose exponents satisfy the equality in the transference inequalities. In this case, we say that the inequalities are \emph{best possible}. A stronger result is to prove that given $k$ exponents $e_{1}, \ldots , e_{k}$, the transference inequalities between these $k$ exponents define a subset of $\mathbb{R}^{k}$ that is exactly the set of all $k$-uples $ (e_{1}(\boldsymbol{\theta}) , \ldots , e_{k}(\boldsymbol{\theta}) )$ as $\boldsymbol{\theta}$ runs through all points $\boldsymbol{\theta}=(\theta_1, \ldots , \theta_n) \in \mathbb{R}^n$ such that $1,\theta_1, \ldots , \theta_n$ are $\mathbb{Q}$-linearly independent. The latter set is called the \emph{spectrum} of the exponents $(e_{1}, \ldots , e_{k})$. \\

 When the dimension is $n=1$, we have the equality $\hat{\omega}_0(\boldsymbol{\theta})=\hat{\omega}(\boldsymbol{\theta})=\hat{\lambda}(\boldsymbol{\theta})=1$. In \cite{JAR}, V. Jarn\'ik showed that in dimension $n=2$, we have the following algebraic relation between $\hat{\omega}_{1}(\boldsymbol{\theta})$ and $\hat{\omega}_0(\boldsymbol{\theta})$:

 \begin{equation}\label{JR} \tag{$\ast$} \hat{\omega}_0(\boldsymbol{\theta})  + \cfrac{1}{ \hat{\omega}_1(\boldsymbol{\theta})} =1 .  \end{equation}
   
   Furthermore, V. Jarn\'ik noted that, in higher dimension $n\geq3$, no algebraic relation holds anymore. He proved \cite[Satz 3]{JAR} that for $n\geq2$, there exist two $n$-tuples of real numbers $\boldsymbol{\theta}=(\theta_1, \ldots \theta_n)$ and $\boldsymbol{\nu}=(\nu_1, \ldots , \nu_n)$ such that

\[  \hat{\omega}_{n-1}(\boldsymbol{\theta}) =\hat{\omega}_{n-1}(\boldsymbol{\nu}) = +\infty ,  \quad   \hat{\omega}_{0}(\boldsymbol{\theta})  = 1 \quad  \textrm{ and }  \quad \hat{\omega}_{0}(\boldsymbol{\nu}) = \cfrac{1}{n-1}.\]

V. Jarn\'ik also proved the following transference theorem:

\begin{thm}[Jarn\'ik, 1938]\label{JarN}
Let $n\geq2$. For any $n$-tuples of real number $\boldsymbol{\theta}=(\theta_1, \ldots \theta_n)$ such that $1,\theta_1, \ldots , \theta_n$ are $\mathbb{Q}$-linearly independent, we have
\[ \cfrac{\hat{\omega}_{n-1}(\boldsymbol{\theta})}{(n-1)\hat{\omega}_{n-1}(\boldsymbol{\theta})+n }  \leq \hat{\omega}_0(\boldsymbol{\theta})  \leq \cfrac{\hat{\omega}_{n-1}(\boldsymbol{\theta})-n+1}{n}. \]
\end{thm}

If $\hat{\omega}_{n-1}(\boldsymbol{\theta})=n$, the interval reduces to the single point $\hat{\omega}_0(\boldsymbol{\theta})=\cfrac{1}{n}$.\\

\begin{rmq}
O. German \cite{OGJar} and A. Khinchin \cite{KhinRmq} claim that  V. Jarn\'ik \cite{JAR} proved the existence of $n$-tuples $\boldsymbol{\theta}=(\theta_1, \ldots , \theta_n)$ with $\hat{\omega}_{n-1}(\boldsymbol{\theta}) = +\infty$ and $\hat{\omega}_0(\boldsymbol{\theta})$ anywhere in the interval $\left[1/(n-1),1\right]$. It appears to the author that this is not written explicitly in \cite{JAR}.\end{rmq}

Recently, O. German \cite{OGJar} improved Theorem \ref{JarN}:
  
  \begin{thm}[German, 2012]\label{OG}
  With the notation of Theorem \ref{JarN}, we have
  \begin{equation} \label{Ger12}\tag{$\ast\ast$}  \cfrac{\hat{\omega}_{n-1}(\boldsymbol{\theta})-1}{(n-1)\hat{\omega}_{n-1}(\boldsymbol{\theta})} \leq  \hat{\omega}_0(\boldsymbol{\theta}) \leq  \cfrac{\hat{\omega}_{n-1}(\boldsymbol{\theta})-(n-1)}{\hat{\omega}_{n-1}(\boldsymbol{\theta})}.\end{equation}
  \end{thm}
  
 Note that the interval reduces to a single point if $n=2$, and that in this case we recover Jarn\'ik's relation \eqref{JR}.\\
 
 The first goal of this paper is to prove that German's inequalities describe the \emph{spectrum} of the two exponents $(\hat{\omega}_{0},\hat{\omega}_{n-1})$ .
   
   \begin{thm}\label{ThQ1}
   Let $n\geq 2$ be an integer, let $\hat{\omega}\in[n,+\infty]$ and let \[\hat{\lambda}\in \left[ \cfrac{\hat{\omega}-1}{(n-1)\hat{\omega}} ,  \cfrac{\hat{\omega}-n+1}{\hat{\omega}} \right],\] where we understand that the interval for $\hat{\lambda}$ is $\left[1/(n-1),1\right]$ when $\hat{\omega}=+\infty$. Then there exist uncountably many $n$-tuples of real numbers $\boldsymbol{\theta}=(\theta_{1}, \ldots, \theta_{n})$,  with $1, \theta_{1}, \ldots, \theta_{n}$ $\mathbb{Q}$-linearly independent,  such that $\hat{\omega}_{n-1}(\boldsymbol{\theta}) = \hat{\omega}$ and $\hat{\omega}_0(\boldsymbol{\theta})=\hat{\lambda}$.  \end{thm}
  
  In \cite{SSJar}, W. Schmidt and L. Summerer obtained independently a similar result, proving that the inequalities \eqref{Ger12} of German are \emph{best possible}.\\

  One can wonder if in higher dimension ($n\geq3$), there exists a Jarn\'ik's-type relation between the $n$ uniform exponents $\hat{\omega}_0, \ldots , \hat{\omega}_{n-1}$. The next theorem states that no such algebraic relation holds.
  
   \begin{thm}\label{ThQ2}
  For every integer $n\geq3$, the $n$ uniform exponents $\hat{\omega}_0, \ldots, \hat{\omega}_{n-1}$ are algebraically independent.   \end{thm}
   
 Thus, the spectrum of the $n$ uniform exponents $\hat{\omega}_0, \ldots, \hat{\omega}_{n-1}$ is a subset of $\mathbb{R}^{n}$ with nonempty interior. \\

  We also know the spectrum of other families of exponents. M. Laurent \cite{ML} described the spectrum of the four exponents $\omega_0,\hat{\omega}_0,\omega_{n-1},\hat{\omega}_{n-1}$ in dimension $n=2$. 
In his PhD thesis, the author gives an alternative proof of this result. However, for $n\geq3$ this spectrum is still unknown. \\

D. Roy showed in \cite{RoySpec} that the going-up and going-down transference inequalities of M. Laurent \cite{MLwd} describe the spectrum of the $n$ exponents $\omega_0, \ldots , \omega_{n-1}$.\\  

 In section \ref{Param}, we introduce Parametric Geometry of Numbers, which is the main tool to prove Theorem \ref{ThQ1}  (section \ref{Dem1}) and Theorem \ref{ThQ2}  (section \ref{dem}), and to give an alternative proof of Theorem \ref{OG} (section \ref{Dem2}) .\\
  

 \section{Parametric Geometry of Numbers}\label{Param}
 
 The Parametric Geometry of Numbers answers a question of W. M. Schmidt \cite{SchLum}. Given a convex body and a lattice, we deform either of them with a one parameter diagonal map. We study the behavior of the successive minima in terms of this parameter. It was developed by W. M. Schmidt and L. Summerer \cite{SchSu,SchSu2}, and further by D.Roy \cite{RoyParam,RoySpec}. Independently, I. Cheung \cite{Cheung1,Cheung2} also developed a similar theory.\\
 
In this paper, we use the notation introduced by D. Roy in \cite{RoyParam,RoySpec} which is essentially dual to the one of W. M. Schmidt and L. Summerer \cite{SchSu,SchSu2}. We refer the reader to these papers for further details. Here $ \boldsymbol{x}\cdot \boldsymbol{y}=x_1y_1+ \cdots + x_ny_n$  is the usual scalar product of vectors $\boldsymbol{x}$ and $\boldsymbol{y}$, and $\|\boldsymbol{x}\|_2=\sqrt{\boldsymbol{x}\cdot \boldsymbol{x}}$ is the usual Euclidean norm.\\

Let $\boldsymbol{u} = (u_0, \ldots , u_n)$ be a vector in $\mathbb{R}^{n+1}$, with Euclidean norm $\|\boldsymbol{u}\|_2=1$. For a real parameter $Q\geq1$ we consider the convex body
\[ \mathcal{C}_{\boldsymbol{u}}(Q) = \left\{ \boldsymbol{x} \in \mathbb{R}^{n+1} \mid \|\boldsymbol{x}\|_2 \leq 1 , \; |\boldsymbol{x} \cdot \boldsymbol{u} | \leq Q^{-1}     \right\}. \]
 For $1\leq d \leq n+1$ we denote by $\lambda_d\left( \mathcal{C}_{\boldsymbol{u}}(Q) \right)$ the $d$-th minimum of $\mathcal{C}_{\boldsymbol{u}}(Q)$ relatively to the lattice $\mathbb{Z}^{n+1}$. For $q\geq0$ and $1\leq d \leq n+1$ we set
\[ L_{\boldsymbol{u},d}(q) = \log   \lambda_d\left( \mathcal{C}_{\boldsymbol{u}}(e^q) \right). \]
Finally, we define the following map associated with $\boldsymbol{u}$:

\[ \begin{array}{rccl}
 \boldsymbol{L_u} :& [0,\infty ) & \to& \mathbb{R}^{n+1}  \\
  & q & \mapsto& (L_{\boldsymbol{u},1}(q), \ldots , L_{\boldsymbol{u},n+1}(q))  .
  \end{array}\]
  The lattice $\mathbb{Z}^{n+1}$ is invariant under permutation of coordinates. Hence,  $\boldsymbol{L_u}$ remains the same if we permute the coordinates in $\boldsymbol{u}$. Since $\|\boldsymbol{u}\|_2=1$ we can thus assume that $u_0\neq0$. \\
  
The following proposition links the exponents of Diophantine approximation associated with $\boldsymbol{\theta}=(\cfrac{u_1}{u_0}, \ldots , \cfrac{u_n}{u_0})$ to the behavior of the map $ \boldsymbol{L_u}$, assuming $u_0\neq0$. It was first stated by W.M. Schmidt and L. Summerer in \cite{SchSu} (Theorem 1.4). It also appears as Relations (1.8) and (1.9) in \cite{SchSu2}. In the notation of D.Roy \cite{RoySpec} (Proposition 3.1), it reads as follows.

\begin{prop}[Schmidt, Summerer, 2009]\label{prop}
Let $\boldsymbol{u} = (u_0, \ldots , u_n) \in \mathbb{R}^{n+1}$, with Euclidean norm $\|\boldsymbol{u}\|_2=1$ and  $u_0\neq 0$. Set $\boldsymbol{\theta}=(\cfrac{u_1}{u_0}, \ldots , \cfrac{u_n}{u_0})$. For $1\leq k \leq n$, we have the following relations:

\begin{eqnarray*}
\liminf_{q\to+\infty}  \cfrac{ L_{\boldsymbol{u},1}(q) + \cdots + L_{\boldsymbol{u},k}(q)}{q} &=& \cfrac{1}{1+{\omega}_{n-k}(\boldsymbol{\theta})},\\
\limsup_{q\to+\infty} \cfrac{ L_{\boldsymbol{u},1}(q) + \cdots + L_{\boldsymbol{u},k}(q)}{q} &=& \cfrac{1}{1+\hat{\omega}_{n-k}(\boldsymbol{\theta})}.\\
\end{eqnarray*}
\end{prop}

Thus, if we know an explicit map $\boldsymbol{P}=(P_1, \ldots, P_{n+1}): [0,\infty) \to \mathbb{R}^{n+1}$, such that 
$\boldsymbol{L}_{\boldsymbol{u}}-\boldsymbol{P}$ is bounded, then we can compute the $2n$ exponents $\hat{\omega}_{0}(\boldsymbol{\theta}), \ldots , \hat{\omega}_{n-1}(\boldsymbol{\theta}), {\omega}_{0}(\boldsymbol{\theta}), \ldots , {\omega}_{n-1}(\boldsymbol{\theta})$ for the above point $\boldsymbol{\theta}$ upon replacing $L_{\boldsymbol{u},i}$ by $P_i$ in the above formulas for $1\leq i \leq n$.\\
For this purpose, we consider the following family of maps, introduced by D. Roy in \cite{RoySpec}.

\begin{defn}[Roy, 2014]
Let $I$ be a subinterval of $[0,\infty)$ with non-empty interior. A generalized $(n+1)$-system on $I$ is a continuous piecewise linear map $\boldsymbol{P} = (P_1, \ldots , P_{n+1}): I \to \mathbb{R}^{n+1}$ with the following three properties.

\begin{description}

\item[(S1)]{For each $q\in I$, we have $0\leq P_1(q) \leq \cdots \leq P_{n+1}(q) $  and $P_1(q) + \cdots + P_{n+1}(q) =q$.}

\item[(S2)]{ If $H$ is a non-empty open subinterval of $I$ on which $\boldsymbol{P}$ is differentiable, then there are integers $\underline{r}, \bar{r}$ with $1\leq \underline{r} \leq \bar{r} \leq n+1$ such that $P_{\underline{r}}, P_{\underline{r}+1}, \ldots , P_{\bar{r}}$ coincide on the whole interval $H$ and have slope $1/(\bar{r}-\underline{r}+1)$ while any other component $P_k$ of $\boldsymbol{P}$ is constant on $H$ .}

\item[(S3)]{ If $q$ is an interior point of $I$ at which $\boldsymbol{P}$ is not differentiable, if $\underline{r}, \bar{r}, \underline{s},\bar{s}$ are the integers for which
\[ P_k'(q^-) = \cfrac{1}{\bar{r}-\underline{r}+1}  \quad  (\underline{r} \leq k \leq \bar{r}) \quad \textrm{  and  }  \quad P_k'(q^+) = \cfrac{1}{\bar{s}-\underline{s}+1}  \quad  (\underline{s} \leq k \leq \bar{s}) \; , \]
and if $\underline{r}<\bar{s}$, then we have $P_{\underline{r}}(q) = P_{\underline{r}+1}(q) = \cdots = P_{\bar{s}}(q)$.
 }\\
\end{description}
\end{defn}

Here $P_k'(q^-)$ (resp. $P_k'(q^+)$) denotes the left (resp. right) derivative of $P_k$ at $q$.
The next result combines Theorem 4.2 and Corollary 4.7 of \cite{RoySpec}.

\begin{thm}[ Roy, 2014]\label{DR}
For each non-zero point $\boldsymbol{u} \in \mathbb{R}^{n+1}$, there exists $q_0\geq 0$ and a generalized $(n+1)$-system $\boldsymbol{P}$ on $[q_0,\infty)$ such that $\boldsymbol{L_u} - \boldsymbol{P}$ is bounded on $[q_0,\infty)$. Conversely, for each generalized $(n+1)$-system $\boldsymbol{P}$ on an interval $[q_0,\infty)$ with $q_0\geq0$, there exists a non-zero point $\boldsymbol{u}\in\mathbb{R}^{n+1}$ such that $\boldsymbol{L_u} - \boldsymbol{P}$ is bounded on $[q_0,\infty)$.\\
\end{thm}

In view of the remark following Proposition \ref{prop}, this result reduces the determination of the joint spectrum of Diophantine approximation exponents to a combinatorial study of generalized $(n+1)$-systems.\\

Although the definition of a generalized $(n+1)$-system $\boldsymbol{P}=(P_1, \ldots , P_{n+1})$ may look complicated, it is easy to understand in terms of the \emph{combined graph} of $\boldsymbol{P}$, that is the union of the graphs of $P_1, \ldots ,P_{n+1}$ over the interval of definition $I$ of $\boldsymbol{P}$. We explain this below. \\

A \emph{division point} of $\boldsymbol{P}$ is an endpoint of $I$ contained in $I$ or an interior point of $I$ at which $\boldsymbol{P}$ is not differentiable. Such points form a discrete subset of $I$. Between two consecutive division points $q^\ast < q$ of $I$, the graph of each component of $\boldsymbol{P}$ is a line segment. All these line segments have slope $0$ except for one line segment of positive slope $1/t$ where $t$ is the number of components  of $\boldsymbol{P}$ whose graph over $[q^\ast,q]$ is that line segment. In view of the condition $P_1 \leq P_2 \leq \cdots \leq P_{n+1}$, there must be consecutive components $P_{\underline{r}}, \ldots , P_{\bar{r}}$ of $\boldsymbol{P}$ with $\bar{r} - \underline{r} +1 =t$. If $q$ is also an interior point of $I$ and if $P_{\underline{s}}, \ldots , P_{\bar{s}}$ are the components of $\boldsymbol{P}$ whose graph has positive slope $\cfrac{1}{\bar{s} - \underline{s} +1}$ to the right of $q$, then there are two cases.\\

\begin{enumerate}[label=\arabic*)]
\item{} If $\underline{r} < \bar{s}$, we say that $q$ is an \emph{ordinary division point}. In this case, we have $P_{\underline{r}}(q) = \cdots = P_{\bar{s}}(q)$ according to \textbf{(S3)}. This implies that $\underline{r} \leq \underline{s}$ and $\bar{r} \leq \bar{s}$. Among $P_{\underline{r}}, \ldots , P_{\bar{s}}$, the components $P_j$ with $\underline{s} \leq j \leq \bar{r}$ (if any) change slope from $\cfrac{1}{\bar{r} - \underline{r} +1}$ to $\cfrac{1}{\bar{s} - \underline{s} +1}$. Those with $j \leq \min(\bar{r},\underline{s}-1)$ change slope from $\cfrac{1}{\bar{r} - \underline{r} +1}$ to $0$. The remaining components $P_j$ with $\bar{r}+1 \leq j \leq \underline{s} -1$ (if any) have constant slope $0$ in a neighborhood of $q$. The reader is invited to draw a picture for himself or to look at those in \cite[\S4]{RoySpec}.
\item{Otherwise}, we have $\underline{r} >\bar{s}$ because it cannot happen that $\underline{r}=\bar{s}$ (or $\boldsymbol{P}$ is differentiable at $q$). Then, we say that $q$ is a \emph{switch point}. In this case, we have $P_{\underline{r}}(q) = \cdots = P_{\bar{r}}(q) > P_{\underline{s}}(q)= \cdots = P_{\bar{s}}(q) $ which mean that the end point of the line segment of slope $\cfrac{1}{\bar{r} - \underline{r} +1}$ at the left of $q$ lies above the initial point of the line segment of slope $\cfrac{1}{\bar{s} - \underline{s} +1}$ at the right of $q$.
\end{enumerate}

It can be shown that the combined graph of a generalized $(n+1)$-system $\boldsymbol{P}$ uniquely determines the map $\boldsymbol{P}$ provided that we know the value of $\boldsymbol{P}$ at one point of its interval of definition. An example of this is shown in \cite[\S4]{RoySpec}. We will see two other examples in the sections \ref{Dem1} and \ref{dem}.\\

In \cite{SchSu,SchSu2} W. M. Schmidt and L. Summerer introduce the following exponents for an integer $1\leq d \leq n+1$:
\[\underline{\varphi}_d = \liminf_{q\to\infty} \cfrac{L_{\boldsymbol{u},d}(q)}{q},\]
\[\bar{\varphi}_d = \limsup_{q\to\infty} \cfrac{L_{\boldsymbol{u},d}(q)}{q}.\]

For these exponents, we have the following analogue of Theorem \ref{ThQ2}:
\begin{thm}\label{ThQ3}
For every integer $n\geq3$, the exponents $\bar{\varphi}_1, \ldots , \bar{\varphi}_{n}$ are algebraically independent.
\end{thm}


\section{Proof of Theorem \ref{ThQ1} } \label{Dem1}

In this section, we construct a family of generalized $(n+1)$-systems. Then, via Theorem \ref{DR}, we get a family of $n$-tuples having the requested properties stated in Theorem \ref{ThQ1}. We first treat the case where $\hat{\omega}_{n-1}$ is finite and $n\geq 3$. We will explain later how to adapt the construction if $n=2$ or $\hat{\omega}_{n-1}$ is infinite.\\

First, note that a generalized $(n+1)$-system with all components equal to $q/(n+1)$ provides via Theorem \ref{DR} a point $\boldsymbol{\theta}$ with $\hat{\omega}_{n-1}(\boldsymbol{\theta})=n$ and $\hat{\omega}_{0}(\boldsymbol{\theta})=1/n$. Thus, we can exclude this case in the next construction.\\

Let $q_0$ be a positive real number, fix a real number $\hat{\omega} > n \geq 2$ and set a parameter $a$ with $\tfrac{1}{n-1} \leq a \leq 1$. We define the sequence $(q_{6m})_{m\geq0}$ by:

\begin{equation*}
q_{6m} = (1+a(\hat{\omega}-n)) q_{6(m-1)} , \;  \textrm{ for }  \; m\geq 1.\\
\end{equation*}

Since $\hat{\omega}>n$, the term $q_{6m}$ goes to infinity as $m$ does.\\

\begin{figure}[!h] 
 \begin{center}
 \begin{tikzpicture}[scale=0.45]
 \draw[black, semithick] (29,1)--(7,1) node [left,black] { $P_2$ } ;
  \draw[black, semithick] (17,2)--(7,2) node [left,black] { $P_3= \cdots = P_n$ } ;
  \draw[black, semithick] (33,4)--(7,4) node [left,black] { $ P_{n+1}$ }; \draw[black, semithick] (33,4) node [right,black] {$P_{1}$};
   \draw[black, semithick] (17,10)--(33,10) node [right,black] {$P_{n+1}$};
   \draw[black, semithick] (29,1)--(33,5) node[midway,above,sloped,black] { slope $1$};\draw[black, semithick] (33,5) node [right,black] {$P_{2}$};
   \draw[black, semithick] (17,10)--(7,0) node [left,black] { $P_1$ } node[midway,above,sloped,black, semithick] { slope $1$};
   \draw[black, semithick] (17,2)--(29,6) node[midway,above,sloped,black] { slope $\cfrac{1}{n-2}$ };
   \draw[black, semithick] (29,6)--(33,6) node [right,black] {$P_3= \cdots = P_n$};
\draw[dashed, black] (8,4)--(8,0) node [below] { $q_{6m}$};
\draw[dashed, black] (9,4)--(9,-1) node [below] { $q_{6m+1}$};
\draw[dashed, black] (11,4)--(11,0) node [below] { $q_{6m+2}$};
\draw[dashed, black] (17,10)--(17,0) node [below] { $q_{6m+3}$};
\draw[dashed, black] (23,10)--(23,0) node [below] { $q_{6m+4}$};
\draw[dashed, black] (29,10)--(29,0) node [below] { $q_{6m+5}$};
\draw[dashed, black] (32,10)--(32,0) node [below] { $q_{6(m+1)}$};
 \end{tikzpicture}
 \end{center}
 \caption{Combined graph of $\boldsymbol{P}$ on a generic interval $[q_{6m},q_{6(m+1)}]$ }\label{fig1}
 \end{figure}
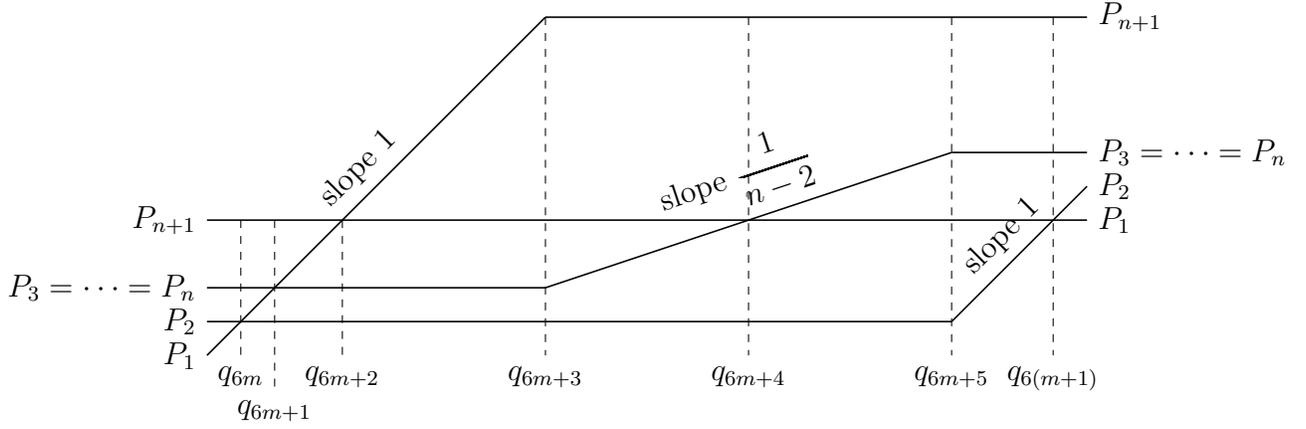

We construct a generalized $(n+1)$-system $\boldsymbol{P}$ whose graph is invariant under the dilation of factor $(1+a(\hat{\omega}-n))>1$ on the interval $[q_0,+\infty)$. Thus, we only need to define $\boldsymbol{P}$ on a generic interval $[q_{6m},q_{6(m+1)}]$. Figure \ref{fig1} shows the pattern of the combined graph of $\boldsymbol{P}$.\\

For every integer $m\geq 0$, we define $\boldsymbol{P}$ at $q_{6m}$ as follows:
\begin{eqnarray*}
P_{1}(q_{6m}) &=&  P_{2}(q_{6m}) = \cfrac{q_{6m}}{\hat{\omega}+1} , \\
P_{3}(q_{6m}) &=& \cdots = P_{n}(q_{6m}) = \cfrac{1 + \left(\cfrac{1-a}{n-2}\right)(\hat{\omega}-n)}{\hat{\omega}+1} q_{6m}, \\
P_{n+1}(q_{6m}) &=& \cfrac{1+a(\hat{\omega}-n)}{\hat{\omega}+1}q_{6m}.
\end{eqnarray*}

Here the parameter $a$ says how large $P_{n+1}$ is at each point $q_{6m}$. The condition $a\geq 1/(n-1)$  imposes the condition $P_{n+1}(q_{6m}) \geq P_{n}(q_{6m})$, and the condition $a\leq1$ imposes that $P_{3}(q_{6m}) \geq P_{2}(q_{6m})$. We have the dilation condition $\boldsymbol{P}(q_{6(m+1)}) = \boldsymbol{P}((1+a(\hat{\omega}-n))q_{6m}) = (1+a(\hat{\omega}-n))\boldsymbol{P}(q_{6m})$ by the definition of the sequence $(q_{6m})_{m\geq0}$.\\

For $k= 0, \ldots , 5$ the graph has only one line segment of positive slope on the interval $[q_{6m+k}, q_{6m+k+1}]$. The graph is clearly the combined graph of a generalized $(n+1)$-system with seven division points $q_{6m}, \ldots , q_{6m+6}$. The points $q_{6m+3}$ and $q_{6m+5}$ are switch points while the others are ordinary division points. Furthermore it is uniquely defined since we know the value of $\boldsymbol{P}$ at the point $q_{6m}$, where as requested \[P_{1}(q_{6m}) + \cdots + P_{n+1}(q_{6m}) = q_{6m}.\] 

Easy computation gives

\[\begin{array}{ll}
q_{6m} = (1+a(\hat{\omega}-n)) q_{6(m-1)} ,& q_{6m+1} = \cfrac{(n-2)(\hat{\omega}+1) + (1-a)(\hat{\omega}-n)}{(n-2) (\hat{\omega}+1)}q_{6m},\\[3mm]
q_{6m+2} = \cfrac{(n+1) + (1+a)(\hat{\omega}-n)}{\hat{\omega}+1}q_{6m} ,& q_{6m+3} = \cfrac{\hat{\omega} + (1+a(\hat{\omega}-n))^{2} }{\hat{\omega}+1}q_{6m},\\[3mm]
q_{6m+4} = \cfrac{1+ (1+a(\hat{\omega}-n))(n+a(\hat{\omega}-n)}{\hat{\omega}+1}q_{6m},& q_{6m+5} = \cfrac{1+2a(\hat{\omega}-n) + \hat{\omega}(1+a(\hat{\omega}-n))}{\hat{\omega}+1}q_{6m}.\\
\end{array}\]\\

We now compute its associated exponents with Proposition \ref{prop}. One can notice that the local extrema of the functions $q\to q^{-1}P_k(q)$, $1\leq k \leq n+1$ are located at division points where $P_k$ changes slope.\\

 Since $\boldsymbol{P}$ is invariant under dilation of factor $C=(1+a(\hat{\omega}-n))$ we have for every $m\geq0$, every $1\leq k \leq n+1$, and every $q$ in $[q_{6m} , q_{6m+6})$ the relation
 
 \[ q^{-1}P_{k}(q) = q^{-1}C^{m}P_{k}\left(q C^{-m}\right),\]
 where $C^{-m}q$ lies in the fundamental interval $[q_{0},q_{6}]$.\\

 Thus, 
 \begin{eqnarray*}
 \limsup_{q\to+\infty} \cfrac{P_{1}(q) }{q} &=& \max_{q_0\leq q \leq  q_6} \cfrac{P_{1}(q) }{q} = \cfrac{P_{1}(q_0) }{q_0} = \cfrac{1}{\hat{\omega}+1},  \\
 \liminf_{q\to+\infty} \cfrac{P_{n+1}(q) }{q} &=& \min_{q_0\leq q \leq  q_6} \cfrac{P_{n+1}(q) }{q} = \cfrac{P_{n+1}(q_2) }{q_2} = \cfrac{1+a(\hat{\omega}-n)}{n+1 +(1+a)(\hat{\omega}-n)} ,\\
 \end{eqnarray*}
 because the component $P_{n+1}$ changes slope from zero to some positive value only at $q_{6m+2}$.\\

Then, according to Proposition \ref{prop}, Theorem \ref{DR} provides an $n$-tuple $\boldsymbol{\theta}=(\theta_{1}, \ldots , \theta_{n})$ such that
\begin{eqnarray*}
\cfrac{1}{\hat{\omega}_{n-1}(\boldsymbol{\theta}) +1} & = &  \limsup_{q\to+\infty} \cfrac{P_{1}(q) }{q} = \cfrac{1}{\hat{\omega}+1},  \\
\cfrac{ \hat{\omega}_0(\boldsymbol{\theta}) }{ \hat{\omega}_0(\boldsymbol{\theta}) +1 } &=& \liminf_{q\to+\infty} \cfrac{P_{n+1}(q) }{q} = \cfrac{1+a(\hat{\omega}-n)}{n+1 +(1+a)(\hat{\omega}-n)} .
\end{eqnarray*}

Thus, this $\boldsymbol{\theta}$ satisfies 

\[  \hat{\omega}_{n-1}(\boldsymbol{\theta}) = \hat{\omega} \quad  \textrm{   and   }  \quad  \hat{\omega}_0(\boldsymbol{\theta}) = \cfrac{1+ a(\hat{\omega}-n)}{\hat{\omega} }. \]

When $a$ runs through the interval $[1/(n-1),1]$, then $\hat{\omega}_0(\boldsymbol{\theta})$ runs through the interval \[\left[ \cfrac{\hat{\omega}-1}{(n-1)\hat{\omega}} ,  \cfrac{\hat{\omega}-(n-1)}{\hat{\omega} } \right].\]

If $n=2$, we remove the line $P_3= \cdots = P_n$ and the interval $[q_{6m+3}, q_{6m+5}]$ from the generic graph on the interval $[q_{6m}, q_{6(m+1)}]$, the parameter $a$ is then forced to be equal to $1$. Thus, we construct $\boldsymbol{\theta}$ with

\[  \hat{\omega}_{1}(\boldsymbol{\theta}) = \hat{\omega} \quad  \textrm{   and   }  \quad  \hat{\omega}_0(\boldsymbol{\theta}) =1-\cfrac{1}{\hat{\omega}},\]
which agrees with Jarn\'ik's relation \eqref{JR}.\\

 If $\hat{\omega}$ is infinite, we replace $\hat{\omega}$ by $m+n+1$ in our construction. For a given real number $q_0$ we consider the sequence $(q_{6m})_{m\geq1}$ defined by
 \[ q_{6m}= (m+1)q_{6(m-1)}.\]
Figure \ref{fig1} still represents the combined graph from $\boldsymbol{P}$ on a generic interval $[q_{6m}, q_{6m+6}]$, with the following settings at $q_{6m}$:
 \begin{eqnarray*}
P_{1}(q_{6m}) &=&  P_{2}(q_{6m}) = \cfrac{q_{6m}}{m+n+2} , \\
P_{3}(q_{6m}) &=& \cdots = P_{n}(q_{6m}) = \cfrac{1 + \left(\cfrac{1-a}{n-2}\right)(m+1)}{m+n+2} q_{6m}, \\
P_{n+1}(q_{6m}) &=& \cfrac{1+a(m+1)}{m+n+2}q_{6m}.
\end{eqnarray*}
 Note that the combined graph is not invariant under dilation anymore. We have
 \begin{eqnarray*}
 \limsup_{q\to+\infty} \cfrac{P_{1}(q) }{q} &=& \limsup_{m\to+\infty} \max_{q_{6m}\leq q \leq q_{6(m+1)}} \cfrac{P_{1}(q) }{q} = \limsup_{m\to+\infty}  \cfrac{P_{1}(q_{6m}) }{q_{6m}} = \limsup_{m\to+\infty}  \cfrac{1}{m+n+2} =0,  \\
 \liminf_{q\to+\infty} \cfrac{P_{n+1}(q) }{q} &=&  \liminf_{m\to+\infty} \min_{q_{6m}\leq q \leq q_{6(m+1)}} \cfrac{P_{n+1}(q) }{q} =  \liminf_{m\to+\infty} \cfrac{P_{n+1}(q_{6m+2} )}{q_{6m+2}} \\ &=& \liminf_{m\to+\infty} \cfrac{1+a(m+1)}{n+1 +(1+a)(m+1)} = \cfrac{a}{a+1} .\\
 \end{eqnarray*}
 
 Again, Theorem \ref{DR} provides us with an $n$-tuple $\boldsymbol{\theta}=(\theta_{1}, \ldots , \theta_{n})$ such that
 
 \[ \hat{\omega}_{n-1}(\boldsymbol{\theta}) = + \infty \quad  \textrm{   and   }  \quad  \hat{\omega}_0(\boldsymbol{\theta}) = a, \]
 where $a$ runs through the interval $[1/(n-1),1]$.\\
 
 Note that if $1,\theta_1, \ldots, \theta_n$ are $\mathbb{Q}$-linearly dependent, then there exists an integer point $\boldsymbol{x}\in\mathbb{Z}^n$ such that $|\boldsymbol{x}\cdot \boldsymbol{u}|=0$. This implies that $L_{\boldsymbol{u},1}(q)$ is bounded above by $\log(\|x\|_2)$. In our construction by dilatation $P_1$ is not bounded, hence the independence by contradiction.\\

 To complete the proof of Theorem \ref{ThQ1}, we have to check that we can construct uncountably many $n$-tuples with given exponents. Let $\hat{\omega}$ and $\hat{\lambda}$ as in Theorem \ref{ThQ1}, and $a$ the parameter such that Theorem \ref{DR} provides an $n$-tuple $\boldsymbol{\theta}$ whose exponents satisfy
 \[\hat{\omega}_{n-1}(\boldsymbol{\theta})= \hat{\omega} \quad \textrm{ and } \quad \hat{\omega}_0(\boldsymbol{\theta})= \hat{\lambda}= \cfrac{1+a(\hat{\omega}-n)}{\hat{\omega}}.\]
 Fix $q_0$ a real number to start the construction from $\boldsymbol{P}$ as above with parameter $a$. For every $\rho_1$ and $\rho_2$ such that $q_0 \leq \rho_1 < \rho_2 \leq q_5$, we denote by $\boldsymbol{P}_{\rho_1}$ and $\boldsymbol{P}_{\rho_2}$ the $(n+1)$-generalized system with parameter $a$ starting in $\rho_1$ and $\rho_2$. We have $\boldsymbol{P}_{\rho_1}(q_6) \neq \boldsymbol{P}_{\rho_2}(q_6)$ and 
 \[  \|\boldsymbol{P}_{\rho_1}({q}_{6m}) - \boldsymbol{P}_{\rho_2}( {q}_{6m}) \|_\infty = \cfrac{q_{6m}}{q_6} \|\boldsymbol{P}_{\rho_1}({q}_{6}) - \boldsymbol{P}_{\rho_2}( {q}_{6}) \|_\infty \to_{n\to\infty} \infty, \]
where $\| (x_1, \ldots , x_n ) \|_\infty = \max_{1\leq k \leq n} |x_k|$. \\

Thus, their difference is unbounded, and  they cannot correspond to the same $\boldsymbol{\theta}$ via Theorem \ref{DR}. \qed


\section{An alternative proof of Theorem \ref{OG}}\label{Dem2}

In this section, we give an alternative proof of Theorem \ref{OG} using arguments from Parametric Geometry of Numbers. As in previous section, we reduce the study of Diophantine properties of a $n$-tuples of real numbers $\boldsymbol{\theta}$ to the study of generalized $(n+1)$-systems. If $\boldsymbol{\theta}=(\theta_1,\dots,\theta_n)\in\mathbb{R}^n$ is such that $1,\theta_1,\dots,\theta_n$ are linearly independent over $\mathbb{Q}$, by Theorem \ref{DR} there exist $q_0>0$ and  a generalized $(n+1)$-system $\boldsymbol{P}=(P_1,\dots,P_{n+1})$ on $[q_0,\infty)$ such that $\boldsymbol{P}-\boldsymbol{L}_{\boldsymbol{u}}$ is bounded where $\boldsymbol{u}=(1,\theta_1,\ldots,\theta_n)$.  Since $\boldsymbol{u}$ has linearly independent coordinates, the first component $P_1$ of $\boldsymbol{P}$ is unbounded.
 For simplicity, we set $\hat{\omega}=\hat{\omega}_{n-1}(\boldsymbol{\theta})$ and $\hat{\lambda}=\hat{\omega}_{0}(\boldsymbol{\theta})$.  Then according to Proposition \ref{prop}, we have
  \begin{equation}\label{formules}
 \limsup_{q\to+\infty} \cfrac{P_{1}(q)}{q} = \cfrac{1}{\hat{\omega}+1} \quad \textrm{ and } \quad   \liminf_{q\to+\infty} \cfrac{P_{n+1}(q)}{q} = \cfrac{\hat{\lambda}}{\hat{\lambda}+1},
 \end{equation}
where we understand that, if $\hat{\omega}=+\infty$, then the limsup  is zero.\\

One can notice that the extremal values of the components of $\boldsymbol{P}$ are reached at the division points. The condition \textbf{(S3)} translates into the fact that for every division point $q$, the right endpoint of the segment with non-zero slope ending at $q$ lies above the left endpoint of the one starting at $q$. A first consequence is that when $P_{1}$ is non constant, it increases until reaching $P_{2}(q)$. A second consequence is the following proposition.\\
 
 \begin{prop}\label{pente}
 For every $1\leq k < m \leq n+1$, if $p_0$ is a point such that $P_k'(p_0^+)>0$, then for every $p>p_0$
 \[ P_m(p) \leq \max\left( P_m(p_0), P_k(p_0) +p-p_0 \right) . \]
 \end{prop}
  In particular, $P_m$ is constant on the interval $[p_0, p_0 + P_m(p_0)-P_k(p_0)]$.\\
  
 The reason is that, if $p_1$ is the largest real number such that $P_m$ is
constant on $[p_0,p_1]$, then the combined graph of $\boldsymbol{P}$ contains a
polygonal line joining the points $(p_0,P_k(p_0))$ and $(p_1,P_m(p_1))$.
Since the line segments composing such a polygonal line have slope in
$[0,1]$, we must have $p_1\geq p_0+P_m(p_0)-P_k(p_0)$.   The conclusion
follows since $P_m(p)\leq \max\{P_m(p_0), P_m(p_0)+p-p_1\}$ for any $p\geq
p_0$.  This is illustrated on the picture below.
  
 
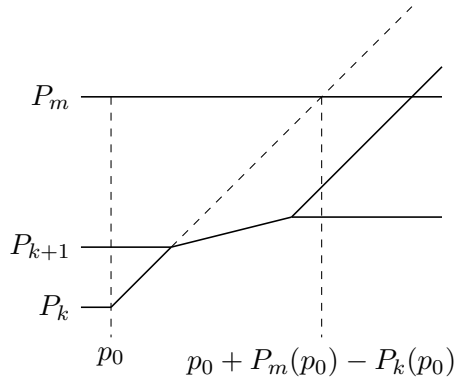
\begin{center}
 \begin{tikzpicture}[scale=0.4]
   \draw[black, semithick] (2,2)--(0,0) ;
   \draw[black, semithick] (2,2)--(-1,2) node [left,black] {$P_{k+1}$};
  \draw[black, semithick] (2,2)--(6,3);
  \draw[black,dashed] (2,2)--(10,10);
  \draw[black, semithick] (7,7)--(-1,7) node [left,black] {$P_m$};
  \draw[black, semithick] (6,3)--(11,3);
     \draw[black, semithick] (6,3)--(11,8);
  \draw[black, semithick] (7,7)--(11,7);
  \draw[black, semithick] (0,0)--(-1,0)node [left,black] {$P_k$};
\draw[dashed, black] (7,7)--(7,-1) node [below] { $p_0 + P_m(p_0)-P_k(p_0)$};
\draw[dashed, black] (0,7)--(0,-1) node [below] { $p_0$};
 \end{tikzpicture}
 \end{center}

The proof of the upper and lower bounds in Theorem \ref{OG} are in the same spirit. For a suitable point of the combined graph, we consider the cases of a polygonal line growing either as fast or as slowly as possible.

\paragraph{Upper bound:}  Suppose first that $\hat{\omega}$ is finite. Let $\varepsilon>0$.  By
\eqref{formules}, there exist arbitrarily large division points $p_0$ where $q^{-1}P_1(q)$ has a local maximum and  \[  \cfrac{1-\varepsilon}{\hat{\omega}+1} \leq \cfrac{ P_{1}(p_0)}{p_0}  \leq   \cfrac{1+\varepsilon}{\hat{\omega}+1}.\]
  
  Since $p_0$ is a local maximum, we have $P_1(p_0)=P_2(p_0)$. Furthermore, $P_{1}(q) \leq P_{2}(q) \leq \cdots \leq P_{n+1}(q)$ and $P_{1}(q) + \cdots + P_{n+1}(q) =q$ provide  
  \[   P_{n+1}(p_0) \leq p_0 -nP_{1}(p_0) \leq   \cfrac{\hat{\omega} +1-n-n\varepsilon}{\hat{\omega}+1}p_0. \]
    At the point $p= p_0 + \cfrac{\hat{\omega}-n-n\varepsilon}{\hat{\omega}+1} p_0$, according to Proposition \ref{pente}, we have the upper bound
      \[P_{n+1}(p) \leq \max(P_{n+1}(p_0),    P_{1}(p_0) + p-p_0) \leq \cfrac{1+\varepsilon + \hat{\omega}-n-n\varepsilon}{\hat{\omega}+1}p_0.\]
Note that equality case corresponds to $\boldsymbol{P}$ with a polygonal line of maximal slope $1$ joining the points $\left(p_{0},P_{1}(p_{0})\right)$ and $\left(p,P_{n+1}(p)\right)$. We deduce that
\[  \cfrac{P_{n+1}(p)}{p}\leq \cfrac{ \hat{\omega} +1-n-(n-1)\varepsilon }{ 2\hat{\omega}-n+1-n\varepsilon }.\]
Since $p$ can be made arbitrarily large, we conclude that
\[ \cfrac{\hat{\lambda}}{\hat{\lambda}+1} = \liminf_{q\to+\infty} \cfrac{P_{n+1}(q)}{q} \leq  \cfrac{ \hat{\omega} +1-n}{ 2\hat{\omega}-n+1},\]
giving that
\[ \hat{\lambda} \leq \cfrac{\hat{\omega}-(n-1)}{\hat{\omega}}.\]

Suppose now that $\hat{\omega}$ is infinite.  Let $\varepsilon>0$. Since $P_1$ is unbounded, there are arbitrarily large values of $p_0$ at which $P_1(p_0)=P_2(p_0)$. At such a point, we have  $P_2'(p_0)>0$.  If $p_0$ is large enough, by \eqref{formules} we also have
\[0 \le \frac{P_1(p_0)}{p_0}\leq \varepsilon\] 
Then, Proposition \ref{pente} applied at the point $p=(2 - n\varepsilon)p_0$ provides
\[P_{n+1}(p) \leq p_0(1-(n-1)\varepsilon).\]
Thus, we get the upper bound
\[  \cfrac{P_{n+1}(p)}{p}\leq \cfrac{p_0(1-(n-1)\varepsilon)}{p_0(2-n\varepsilon)}.\]
Since $p$ can be made arbitrarily large, we conclude that
\[ \hat{\lambda} \leq 1.\]

Hence, we have proved the upper bound in Theorem \ref{OG}.

\paragraph{Lower bound:} 

If $P_1(q)=P_{n+1}(q)$ for arbitrarily large $q$, then $\hat{\omega}=n$ and $\hat{\lambda}=1/n$, and the inequalities of Theorem \ref{ThQ1} are satisfied. So, we may assume that $P_1(q)<P_{n+1}(q)$ for any sufficiently large $q$.\\

  Suppose first that $\hat{\omega}$ is finite. Let $\varepsilon_{1}>0$. By \eqref{formules}, there exists a real number $q_{0}$ such that $q\geq q_{0}$ implies
  \begin{eqnarray*}
 \cfrac{P_{1}(q)}{q} \leq \cfrac{1+\varepsilon_{1}}{\hat{\omega}+1} \qquad \textrm{  and  } \qquad P_{1}(q)\neq P_{n+1}(q). 
  \end{eqnarray*}
   
  Let $\varepsilon_{2}>0$. There exist arbitrarily large division points $p\geq q_{0}$ where $q^{-1}P_{n-1}(q)$ has a local minimum and 
  \[  \left| \cfrac{P_{n+1}(p)}{p} - \cfrac{\hat{\lambda}}{\hat{\lambda}+1} \right|\leq \varepsilon_{2}.  \]
Let $p_0 = \max\left\{q\leq p \mid P_{1}(q)=P_{2}(q)\right\}$. At the point $p_0$ we have 
\[P_{1}(p_0) = P_{2}(p_0) \leq \cfrac{1+\varepsilon_{1}}{\hat{\omega}+1}p_0 \quad \textrm{ and } \quad P_{n+1}(p_0) \geq  \cfrac{p_0-2P_{1}(p_0)}{n-1},\]
since $p_0= P_1(p_0)+ \cdots + P_{n+1}(p_0) \leq 2P_1(p_0)+ (n-1)P_{n+1}(p_0)$.\\

We first show that $q\to P_{1}(q)$ is constant on the interval $[p_0,p]$. If not, there exists a real number $p_0<p_1<p$ where $P_{1}$ has slope $>0$.  Since $p$ is a local minimum from $q^{-1}P_{n+1}(q)$, then $P_{n+1}$ changes slope at $p$. Then, $P_1(p)\neq P_{n+1}(p)$ and condition \textbf{(S3)} imply that $P_{1}'(p^-)=0$. Thus, there exists a point in the interval $(p_1,p)$ where $P_1$ changes slope from $>0$ to $0$. At this point $P_1=P_2$, which contradicts the definition of $p_0$. Thus,
\[ P_{1}(p_0) = P_{1}(p).\] 
We can write
\begin{eqnarray*}
p =  \sum_{k=1}^{n+1}P_{k}(p) &\leq& nP_{n+1}(p) + P_{1}(p_0).
\end{eqnarray*}
Note that equality provides that all components except $P_{1}$ are equal. In this case, we have a polygonal line joining $(p_{0},P_{1}(p_{0})$ and $(p,P_{n+1}(p))$ growing as slowly as possible.\\

We deduce the lower bound
\[ \cfrac{P_{n+1}(p)}{p} \geq \cfrac{P_{n+1}(p)}{nP_{n+1}(p) + P_{1}(p_0)},  \]
where the right hand side is an increasing function of $P_{n+1}(p)$. Since \[P_{n+1}(p) \geq P_{n+1}(p_0) \geq \cfrac{p_0-2P_{1}(p_0)}{n-1},\] we have
\begin{eqnarray*}
 \cfrac{P_{n+1}(p)}{p} &\geq& \cfrac{p_0 - 2P_{1}(p_0)}{np_0-(n+1)P_{1}(p_0)},\\
\end{eqnarray*}
where the right hand side is a decreasing function of $P_{1}(p_0)$. Since \[P_{1}(p_0) \leq \cfrac{1+\varepsilon_{1}}{\hat{\omega}+1}p_0 ,\] 
we have
\begin{eqnarray*}
\cfrac{P_{n+1}(p)}{p} &\geq&  \cfrac{\hat{\omega} -1 -2\varepsilon_{1} }{n\hat{\omega} - 1 -(n+1)\varepsilon_{1}}. \\
\end{eqnarray*}
Finally,
\[ \cfrac{\hat{\lambda}}{\hat{\lambda}+1}  \geq \cfrac{\hat{\omega} -1 -2\varepsilon_{1} }{n\hat{\omega} - 1 -(n+1)\varepsilon_{1}} -\varepsilon_{2}. \]
This gives the expected bound
\[\hat{\lambda} \geq \cfrac{\hat{\omega}-1}{(n-1)\hat{\omega}}.\]

Suppose now that $\hat{\omega}$ is infinite. Choose $q_0$ so that $q\geq q_{0}$ implies
  \begin{eqnarray*}
 0 \leq \cfrac{P_{1}(q)}{q} \leq \varepsilon_{1} \quad \textrm{ and } \quad P_{1}(q) \neq P_{n+1}(q). 
  \end{eqnarray*}
  Following the same steps as in the finite case, with the same choice of $p$ we obtain :
  \[\cfrac{P_{n+1}(p)}{p} \geq \cfrac{1-2\varepsilon_1}{n-(n+1)\varepsilon_1}.\]
  Thus, we get  \[   \cfrac{\hat{\lambda}}{\hat{\lambda}+1} \geq \cfrac{1-2\varepsilon_1}{n-(n+1)\varepsilon_1}-\varepsilon_2.\]
This gives the expected lower bound \[\hat{\lambda} \geq \cfrac{1}{n-1}.\] \qed



\section{Proof of Theorems \ref{ThQ2} and \ref{ThQ3} }\label{dem}

In this section, we construct a family of generalized $(n+1)$-systems depending on $n$ parameters which via Theorem \ref{DR} provides us with a family of $n$-tuples $\boldsymbol{\theta}$ whose uniform exponents are expressed as a function of these $n$ parameters. Then, we show that these functions are algebraically independent.\\

Fix the dimension $n\geq3$. Choose $n+2$ parameters $A_{1}, A_{2}, \ldots , A_{n+1},C$ satisfying 

\begin{equation}\label{a}\tag{0}
\begin{split}
 0< A_1&= A_2 < A_3 < A_4 < \cdots <A_{n+1} ,\\
1&=A_1+A_2+ \cdots + A_{n+1} ,\nonumber\\
\cfrac{A_{k+1}}{A_k} &< C \; < \cfrac{A_{k+2}}{A_k} \textrm{ for } 2\leq k \leq n-1,\\
1&<\cfrac{A_{n+1}}{A_n} < C.
\end{split}
\end{equation}

\begin{figure}[!h] 
 \begin{center}
 \begin{tikzpicture}[scale=0.5]
 
  \fill (1,1) circle[radius=3pt];
\fill (5,4) circle[radius=3pt];
\fill (13,10) circle[radius=3pt];
\fill (16,12) circle[radius=3pt];
\fill (20,16) circle[radius=3pt];
\fill (22,17) circle[radius=3pt];
\fill (23,18) circle[radius=3pt];

\draw (9,8) node {$\cdots$};
\draw (9,15) node {$\cdots$};
\draw (9,3) node {$\cdots$};

\draw (18,8) node {$\cdots$};
\draw (18,15) node {$\cdots$};
\draw (18,3) node {$\cdots$};
 
\draw[black, semithick] (3,3)--(0,0) node [left,black] {};
\draw[black, semithick] (8,1)--(0,1) node [left,black] {$A_1=A_2$};
\draw[black, semithick] (3,3)--(0,3) node [left,black] {$A_3$};
\draw[black, semithick] (3,3)--(5,4) node [left,black] {};
\draw[dashed,black] (5,4)--(0,4) node [left,black] {$CA_2$};
\draw[black, semithick] (5,4)--(8,4) node [left,black] {};
\draw[black, semithick] (5,4)--(6,5) node [left,black] {};
\draw[black, semithick] (6,5)--(0,5) node [left,black] {$A_4$};
\draw[black, semithick] (8,9)--(0,9) node [left,black] {$A_k $}node [midway,above,black] {$P_{k}$};
\draw[black, semithick] (8,11)--(0,11) node [left,black] {$A_{k+1}$} node [midway,above,black] {$P_{k+1}$};
\draw[black, semithick] (8,16)--(0,16) node [left,black] {$A_{n+1}$}node [midway,above,black] {$P_{n+1}$};
\draw[black, semithick] (6,5)--(8,6) node [left,black] {};

\draw[black, semithick] (10,8)--(11,9) node [left,black] {};
\draw[black, semithick] (11,9)--(13,10) ;
\draw[black, semithick] (11,9)--(10,9) node [left,black] {};
\draw[black, semithick] (13,10)--(17,10) node [below,midway,black] {$P_{k-1}$};
\draw[black, semithick] (13,10)--(14,11) ;

\draw[black, semithick] (14,11)--(10,11) node [above,midway,black] {$P_{k+1}$};
\draw[black, semithick] (14,11)--(16,12) ; 
\draw[black, semithick] (16,12)--(17,12) ;
\draw[black, semithick] (16,12)--(17,13) node [left,black] {};
\draw[black, semithick] (10,1)--(17,1);
\draw[black, semithick] (10,4)--(17,4);
\draw[black, semithick] (10,16)--(17,16);

\draw[black, semithick] (19,15)--(20,16) node [left,black] {};
\draw[black, semithick] (20,16)--(19,16) node [left,black] {};
\draw[black, semithick] (20,16)--(22,17) node [left,black] {};
\draw[black, semithick] (22,17)--(22.5,17);
\draw[black, semithick] (22.5,17)--(27,17) node [midway, above,black] {$P_{n}$};
\draw[black, semithick] (22,17)--(23,18) ;
\draw[black, semithick] (23,18)--(27,18) node [midway,above,black] {$P_{n+1}$};
\draw[dashed, black] (23,18)--(0,18) node [left,black] {$CA_{n+1}$};
\draw[dashed, black] (22,17)--(0,17) node [left,black] {$CA_{n}$};
\draw[black, semithick] (23,1)--(27,5) node [below, midway,black] {$P_1$};
\draw[black, semithick] (19,1)--(23,1);
\draw[black, semithick] (22,4)--(27,4) node [above,midway,black] {$P_2$};
\draw[black, semithick] (19,4)--(22,4);
\draw[black, semithick] (22,10)--(27,10) node [above,midway,black] {$P_{k-1}$};
\draw[black, semithick] (19,10)--(22,10);
\draw[dashed, black] (13,10)--(0,10) node [left,black] {$CA_{k-1}$};
\draw[black, semithick] (19,12)--(22,12);
\draw[black, semithick] (22,12)--(27,12) node [midway,above,black] {$P_{k}$};

\draw[dashed, black] (1,18)--(1,0) node [below] { $1$};
\draw[dashed, black] (26,18)--(26,0) node [below] { $C$};
\draw[dashed, black] (11,9)--(11,0) node [below] { $\delta_{k-1,1}$};
\draw[dashed, black] (13,10)--(13,-1) node [below] { $\delta_{k-1,2}$};
\draw[dashed, black] (14,11)--(14,0) node [below] { $\delta_{k,1}$};
\draw[dashed, black] (16,12)--(16,0) node [below] { $\delta_{k,2}$};
\draw[dashed, black] (3,3)--(3,0) node [below] { $\delta_{2,1}$};
\draw[dashed,black] (5,4)--(5,0) node [below] {$\delta_{2,2}$};
\draw[dashed,black] (6,5)--(6,-1) node [below] {$\delta_{3,1}$};
\draw[dashed, black] (23,18)--(23,0) node [below] { $\delta_{n+1,1}$};

 \end{tikzpicture}
 \end{center}
 \caption{Pattern of the combined graph of $\boldsymbol{P}$ on the fundamental interval $[1,C]$. }\label{fig2}
 \end{figure}
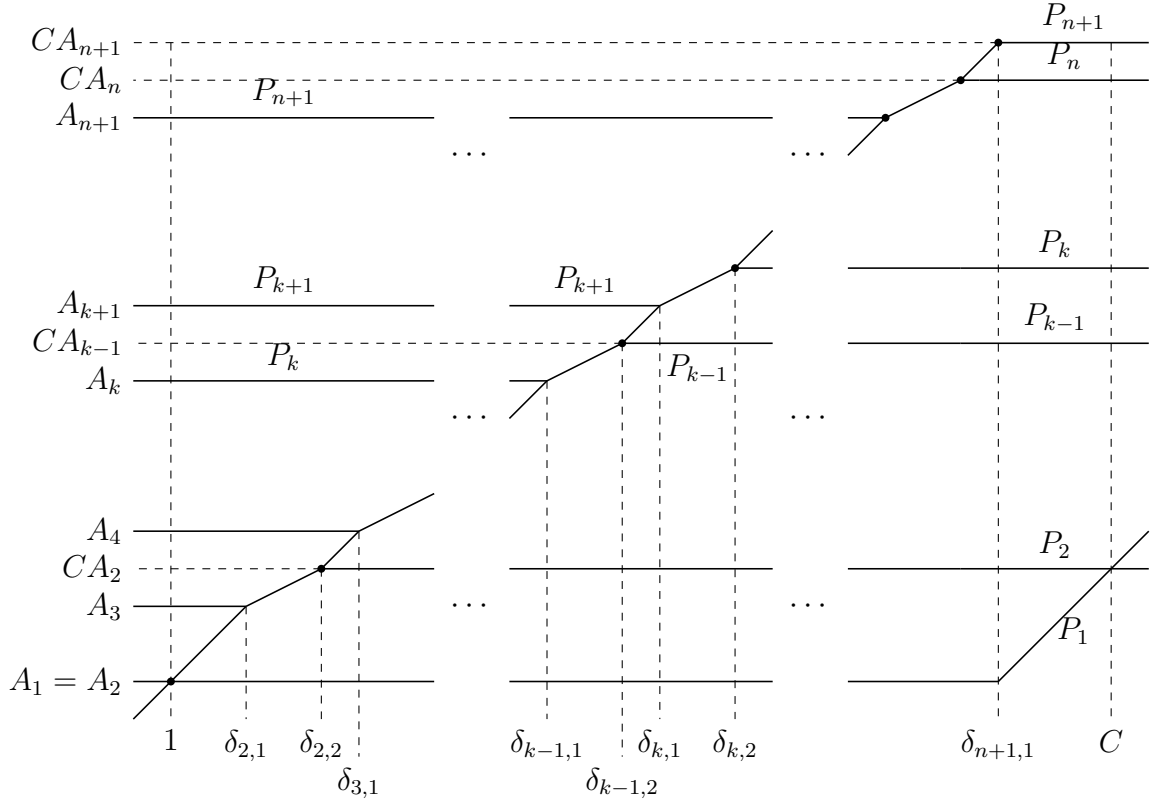

We consider the generalized $(n+1)$-system $\boldsymbol{P}$ on the interval $[1,C]$ whose combined graph is given by Figure \ref{fig2}, where
\[ P_k(1) = A_k \quad \textrm{ and } \quad P_{k}(C)=CA_{k} \quad \textrm{ for } \quad 1\leq k \leq n+1.\]
On each interval between two consecutive division points, there is only one line segment with nonzero slope.
This line segment has slope $1$ on the intervals $[1,\delta_{2,1}]$, $[\delta_{n+1,1},C]$ and $[\delta_{k-1,2},\delta_{k,1}]$ for $3\leq k \leq n+1$, and has slope $1/2$ on the interval $[\delta_{k,1},\delta_{k,2}]$ , for $3\leq k \leq n$ , where the two components $P_{k}$ and $P_{k+1}$ coincide. We have $ 2n+1$ division points $1$, $C$, $\delta_{k,1}$ and $\delta_{l,2}$ for $2\leq k \leq n+1$ and $2 \leq l \leq n$. They are all ordinary division points except $\delta_{n+1,1}$ which is a switch point. Note that the conditions \eqref{a} are consistent with the graph. The points which will be most relevant for the proofs are labeled with black dots. \\

We extend $\boldsymbol{P}$ to the interval $[1,\infty)$ by self-similarity, that is $\boldsymbol{P}(q) = C^{m}\boldsymbol{P}(qC^{-m})$ for every positive integer $m$. In view of the value of $\boldsymbol{P}$ and its derivative at $1$ and $C$, one sees that this extension provides a generalized $(n+1)$-system on $[1,\infty)$.\\

Proposition \ref{prop} suggests to define quantities $\hat{W}_{n-1}, \ldots , \hat{W}_{0}$ by
 \begin{align}\label{W}
 \cfrac{1}{1+\hat{W}_{n-k}}  &:= \limsup_{q\to+\infty} \cfrac{P_{1}(q) + \cdots + P_k(q)}{q} \; , \quad 1\leq k \leq n .
  \end{align}
 Since $\boldsymbol{P}$ is invariant under dilation of factor $C$, we can replace $\limsup_{q\to\infty}$ by $\max_{[1,C]}$ in the above formulae.\\
  
We observe that for $1\leq k \leq n$, the function $P_1 + \cdots + P_k$ has slope $1$ on the intervals $[1,\delta_{k,1}]$ and $[\delta_{n+1,1},C]$, slope $1/2$ on the interval $[\delta_{k,1},\delta_{k,2}]$ and is constant on the interval $[\delta_{k,2},\delta_{n+1,1}]$. Thus the maximum on $[1,C]$ of the function $q\to q^{-1}(P_1(q) + \cdots + P_k(q))$ is reached either at $\delta_{k,1}$ or at $\delta_{k,2}$, when slope changes from $1$ to $1/2$ or  from $1/2$ to $0$. Namely, the maximum is reached at $\delta_{k,1}$ if 
\begin{equation}\label{delta}
\cfrac{P_{1}(\delta_{k,1}) + \cdots +P_k(\delta_{k,1}) }{\delta_{k,1}}\geq\cfrac{1}{2}
\end{equation}
and at $\delta_{k,2}$ if the lefthand side is $\leq 1/2$. We deduce that for $1 \leq k \leq n$, 
\[\hat{W}_{n-k} = \cfrac{P_{k+1}(q) + \cdots + P_{n+1}(q)}{P_{1}(q) + \cdots + P_{k}(q)} \quad \textrm{ where } \quad q=\left\{ \begin{array}{cl} \delta_{k,1} &\textrm{ if \eqref{delta} is satisfied }\\ \delta_{k,2} & \textrm{ otherwise} \end{array}\right. .\]

For $2\leq k \leq n+1$, we have the following values at $\delta_{k,1}$ and $\delta_{k,2}$:

\[\begin{array}{cc}
 P_{i}(\delta_{k,1}) = \left\{\begin{array}{cl} A_{1} & \textrm{ if } i=1\\ CA_{i} & \textrm{ if } 2\leq i \leq k-1 \\ A_{k+1}  &\textrm{ if }   i = k \\ A_{i} & \textrm{ if } k+1 \leq i \leq n+1\end{array}\right. ,& 
 P_{i}(\delta_{k,2}) = \left\{\begin{array}{cl} A_{1} & \textrm{ if } i=1\\ CA_{i} & \textrm{ if } 2\leq i \leq k \\ CA_{k}  &\textrm{ if }   i = k+1 \\ A_{i} & \textrm{ if } k+2 \leq i \leq n+1\end{array}\right. .
\end{array}\]

 It is easy to check that the parameters
  \begin{equation}\label{p}
C=3, \quad A_{1}=A_{2}=2^{-n}, \quad  A_{k}=2^{-n+k-2} \quad \textrm{ for } \; 3 \leq k \leq n+1
 \end{equation}
satisfy the conditions \eqref{a}. For this choice of parameters, the lefthand side of inequality \eqref{delta} is $>1/2$ for $1\leq k \leq n-1$ and  $<1/2$ for $k=n$. This property remains true for $(C,A_{2}, \ldots , A_{n})$ in an open neighborhood of $(3,2^{-n}, \ldots , 2^{-2})$ provided that we set $A_{1}=A_{2}$ and $A_{n+1}=1-(A_{1}+ \cdots + A_{n})$. In this neighborhood, the quantities $\hat{W}_{0}, \ldots , \hat{W}_{n-1}$ are given by the following rational fractions in $\mathbb{Q}(C,A_{2}, A_{3}, \ldots , A_{n})$ :
 \begin{align}
\begin{split}
\hat{W}_{n-1} &= \cfrac{1}{A_{2}} -1, \\ 
\hat{W}_{n-k} &=   \cfrac{1- (2A_2+ A_{3}+ A_{4} + \cdots + A_{k+1})+CA_{k}}{A_2 + C(A_2 + \cdots + A_{k})} \; , \quad 2 \leq k \leq n-1 \\
\hat{W}_{0} & = \cfrac{1-(2A_{2}+A_{3} + A_{4}+ \cdots +A_{n}) }{A_2 + C(A_2 + \cdots + A_{n-1}) }. \\
 \end{split}
\end{align}

Since $\hat{W}_{0}, \ldots , \hat{W}_{n-1}$ come from a generalized $(n+1)$-system $\boldsymbol{P}$, Theorem \ref{DR} provides a point $\boldsymbol{\theta}$ in $\mathbb{R}^{n}$ such that $\hat{\omega}_{k}(\boldsymbol{\theta})= \hat{W}_{k}$ for every $0\leq k \leq n-1$. Thus, to prove Theorem \ref{ThQ2}, it is sufficient to show that the rational fractions $\hat{W}_{0}, \ldots , \hat{W}_{n-1}\in \mathbb{Q}(C,A_{2}, A_{3}, \ldots , A_{n})$ are algebraically independent.\\

Suppose on the contrary that there exists an irreducible polynomial $R\in\mathbb{Q}(X_{1}, \ldots , X_{n})$ such that 
\[ R\left(\hat{W}_{0}, \hat{W}_{1}, \ldots , \hat{W}_{n-1}\right)=0.\]
Specializing $C$ in $0$, we obtain
\[R\left( \cfrac{1-A_{2}-A_{2}- \cdots - A_{n}}{A_{2}}, \cfrac{1-A_{2}-A_{2}- \cdots - A_{n}}{A_{2}}, \ldots , \cfrac{1-A_{2}-A_{2}-A_{3}}{A_{2}}, \cfrac{1-A_{2}}{A_{2}}\right)=0 .\]
Here, the first two rational fractions are the same, and the last $n-1$ rational fractions generate the field $\mathbb{Q}(A_{2}, A_{3}, \ldots , A_{n})$. Therefore the latter are algebraically independent, and $R= \alpha (X_{2}-X_{1})$ for a nonzero constant $\alpha\in\mathbb{Q}$. This is impossible since $\hat{W}_{0}\neq \hat{W}_{1}$.\qed

\newpage
{\it Proof of Theorem \ref{ThQ3}}\\

We consider the same generalized $(n+1)$-system as above. Notice that for $1\leq k \leq n$ we have $P_k\leq P_{n+1}$ and therefore
\[ 0 \leq  \cfrac{P_k(q)}{q} \leq 1/2. \]
Since all nonzero slopes of the combined graph are at least $1/2$, the maxima of the functions $q \mapsto q^{-1}P_{k}(q)$ are reached at points where $P_{k}$ changes slope from $1$ or $1/2$ to $0$. It happens that for each component there is only one such point on the interval $[1,C[$.\\ 
The definition of the exponents $\bar{\varphi}_{k}$ leads to define quantities $F_{k}$ by 

\[F_{k} := \limsup_{q\to\infty} \cfrac{P_k(q)}{q} = \max_{[1,C]} \cfrac{P_k(q)}{q} = \cfrac{P_k(p)}{p} \quad \textrm{ where } \quad p = \left\{ \begin{array}{cl} 1 &  \textrm{ if } k=1,\\ \delta_{k,2} &\textrm{ if } 2\leq k \leq n. \end{array}\right. \]

We express the quantities $F_{1}, \ldots F_{n}$ as rational fractions in $\mathbb{Q}(C,A_{2}, \ldots , A_{n})$, using the relations $A_{1}=A_{2}$ and $A_{n+1}= 1- A_{1}-A_{2}- \cdots - A_{n}$ :

\begin{eqnarray*}
F_{1} &=& A_{1},\\
 F_k &=&  \cfrac{CA_k}{A_1 + C(A_2+ \cdots + A_k) + CA_k + 1- (2A_{2}+ A_{3} + \cdots + A_{k+1})}.\\
 \end{eqnarray*}

Since $F_{1}, \ldots , F_{n}$ come from a generalized $(n+1)$-system $\boldsymbol{P}$, by Theorem \ref{DR} there exists a point $\boldsymbol{\theta}$ in $\mathbb{R}^{n}$ such that $\bar{\varphi}_k(\boldsymbol{\theta})= F_{k}$ for every $1\leq k \leq n$. To prove Theorem \ref{ThQ3} it is sufficient to show that the rational fractions $F_{1}, \ldots , F_{n}\in \mathbb{Q}(C,A_{2}, A_{3}, \ldots , A_{n})$ are algebraically independent.\\

Suppose that there exists an irreducible polynomial $R\in\mathbb{Q}(X_{1}, \ldots , X_{n})$ such that
\[ R(F_{1}, \ldots , F_{n})=0.\]
Specializing $C$ in infinity, we obtain
\[R\left(A_{2}, \cfrac{1}{2}, \cfrac{A_{3}}{(A_{2}+A_{3})+A_{3}}, \ldots , \cfrac{A_{n}}{(A_{2}+ \ldots + A_{n})+A_{n}}\right)=0\]
where all coordinates except $1/2$ are algebraically independent. Thus, $R$ is a constant multiple of $2X_{2}-1$, which contradicts $F_{2}\neq 1/2$.\qed

We are not able to prove Theorem \ref{ThQ3} for the $n+1$ exponents $\bar{\varphi}_1, \ldots , \bar{\varphi}_{n+1}$ with this construction. However with some extra work, we can show that the theorem holds for any $n$ exponents among them.\\

\section*{Acknowledgements}

The author would like to thank the referee and Damien Roy for careful reading and useful remarks to simplify and shorten the proofs.

 \bibliographystyle{plain}
 \bibliography{Biblio}

\begin{thebibliography}{10}

\bibitem{BugLau2}
Yann Bugeaud and Michel Laurent.
\newblock On exponents of homogeneous and inhomogeneous diophantine
  approximation.
\newblock {\em Moscow Math. J.}, 5:747--766, 2005.

\bibitem{BugLau}
Yann Bugeaud and Michel Laurent.
\newblock Exponents of diophantine approximation.
\newblock {\em Diophantine Geometry Proceedings}, 4:101--121, 2007.

\bibitem{Cheung1}
Yitwah Cheung.
\newblock Special divergent trajectories for a homogeneous flow.
\newblock 2008.
\newblock
  \href{http://www.math.uchicago.edu/~geometry/gt_seminar.F2008.html}{Seminar
  of Geometry, Chicago University}.

\bibitem{Cheung2}
Yitwah Cheung.
\newblock Prescriptions for a diagonal flow on the space of lattices.
\newblock 2015.
\newblock
  \href{http://mjcnt.phystech.edu/conference/aarhus/abstracts/Cheung.pdf}{Diophantine
  Approximation and Related Topics, Aarhus}.

\bibitem{OGJar}
Oleg~N. German.
\newblock On {D}iophantine exponents and {K}hintchine's transference principle.
\newblock {\em Mosc. J. Comb. Number Theory}, 2(2):22--51, 2012.

\bibitem{JAR}
Vojt{\v e}ch Jarn{\'\i}k.
\newblock Zum khintchineschen "{{\"U}}bertragungssatz".
\newblock {\em Trav. Inst. Math. Tbilissi}, 3:193--212, 1938.

\bibitem{Khin2}
Alexander~Ya. Khinchin.
\newblock {\"U}ber eine klasse linearer diophantischer approximationen.
\newblock {\em Rend. Circ. Mat. Palermo 50}, pages 170--195, 1926.

\bibitem{Khin1}
Alexander~Ya. Khinchin.
\newblock Zur metrischen theorie der diophantischen approximationen.
\newblock {\em Math.Z.}, 24:706--714, 1926.

\bibitem{KhinRmq}
Alexander~Ya. Khinchin.
\newblock On some applications of the method of the additional variable.
\newblock {\em Amer. Math. Soc. Translation}, 1950(18):14, 1950.

\bibitem{ML}
Michel Laurent.
\newblock Exponents of diophantine approximmation in dimension two.
\newblock {\em Canad. J. Math.}, 61:165--189, 2009.

\bibitem{MLwd}
Michel Laurent.
\newblock On transfer inequalities in {D}iophantine approximation.
\newblock In {\em Analytic number theory}, pages 306--314. Cambridge Univ.
  Press, Cambridge, 2009.

\bibitem{Mink}
Hermann Minkowski.
\newblock {\em Geometrie der {Z}ahlen}.
\newblock Bibliotheca Mathematica Teubneriana, Band 40. Johnson Reprint Corp.,
  New York-London, 1968.

\bibitem{RoyParam}
Damien Roy.
\newblock On {S}chmidt and {S}ummerer parametric geometry of numbers.
\newblock {\em Ann. of Math.}, 182:739--786, 2015.

\bibitem{RoySpec}
Damien Roy.
\newblock Spectrum of the exponents of best rational approximation.
\newblock {\em Math. Z.}, 283:143--155, 2016.

\bibitem{SchH}
Wolfgang~M. Schmidt.
\newblock On heights of algebraic subspaces and diophantine approximations.
\newblock {\em Ann. of Math. (2)}, 85:430--472, 1967.

\bibitem{SchLum}
Wolfgang~M. Schmidt.
\newblock Open problems in {D}iophantine approximation.
\newblock In {\em Diophantine approximations and transcendental numbers
  ({L}uminy, 1982)}, volume~31 of {\em Progr. Math.}, pages 271--287.
  Birkh\"auser Boston, Boston, MA, 1983.

\bibitem{SchSu}
Wolfgang~M. Schmidt and Leonhard Summerer.
\newblock Parametric geometry of numbers and applications.
\newblock {\em Acta Arithmetica}, 140(1):67--91, 2009.

\bibitem{SchSu2}
Wolfgang~M. Schmidt and Leonhard Summerer.
\newblock Diophantine approximation and parametric geometry of numbers.
\newblock {\em Monatsch. Math.}, 169(1):51--104, 2013.

\bibitem{SSJar}
Wolfgang~M. Schmidt and Leonhard Summerer.
\newblock The generalization of jarnik's identity.
\newblock {\em Acta Arithmetica}, to appear.

\end{thebibliography}

 \end{document}